\documentclass[12pt]{article}

\usepackage{amsmath,amssymb,url}
\usepackage{tikz,color}

\newtheorem{theorem}{Theorem}[section]
\newtheorem{corollary}[theorem]{Corollary}
\newtheorem{lemma}[theorem]{Lemma}
\newtheorem{proposition}[theorem]{Proposition}

\newtheorem{problem}[theorem]{Problem}

\newtheorem{remark}[theorem]{Remark}

\newcommand{\tr}{{\rm Tr}}

\newcommand{\proof}{\noindent{\bf Proof.\ }}
\newcommand{\qed}{\hfill $\square$ \bigskip}
\newcommand{\smallqed}{{\tiny ($\Box$)}}

\vspace{8mm}
\title{Constructing  new families of transmission irregular graphs}

\author{Kexiang Xu$\/^{a}$, Sandi Klav\v{z}ar$^{b}$ \\\\
 $^{a}$ \small  College of Science, Nanjing University of
 Aeronautics \& Astronautics,\\
 \small Nanjing, Jiangsu 210016, PR China\\
$^{b}$ \small Faculty of Mathematics and Physics, University of Ljubljana, Slovenia\\
\small {\tt kexxu1221@126.com} (K.\ Xu)\\
\small {\tt sandi.klavzar@fmf.uni-lj.si} (S.\ Klav\v{z}ar)
}

\begin{document}

\maketitle

\begin{abstract}
The transmission of a vertex $v$ of a graph $G$ is the sum of distances from $v$ to all the other vertices in $G$. A graph  is transmission irregular if all of its vertices have pairwise different transmissions.   A starlike tree $T(k_1,\ldots,k_t)$ is a tree obtained by attaching to an isolated vertex $t$ pendant paths of lengths $k_1,\ldots,k_t$, respectively. It is proved that if a starlike tree $T(a,a+1,\ldots,a+k)$, $k\ge 2$, is of odd order, then it is transmission irregular. $T(1,2,\ldots,\ell)$, $\ell \ge 3$, is transmission irregular if and only if $\ell \notin \{r^2 + 1:\ r\ge 2\}$. Additional infinite families among the starlike trees and bi-starlike trees are determined. Transmission irregular unicyclic infinite families are also presented, in particular, the line graph of $T(a,a+1,a+2)$, $a\ge 2$, is transmission irregular if and only if $a$ is even.
\end{abstract}


\medskip\noindent
\textbf{AMS Math.\ Subj.\ Class.\ (2010)}: 05C12, 05C76

\section{Introduction}

If $G = (V(G), E(G))$ is a graph, we use the notations $n(G)=|V(G)|$ and $m(G)=|E(G)|$, and denote by $d_G(u,v)$ the shortest-path distance between vertices $u, v\in V(G)$. The {\em transmission} $\tr_G(v)$ (or \tr(v) for short if the graph $G$ is clear from the context) of a vertex $v\in V(G)$ is the sum of distances from $v$ to the vertices in $G$, that is,
$$\tr_G(v)=\sum\limits_{u\in V(G)}d_G(u,v)\,.$$
With this notation we have $W(G)=\frac{1}{2}\sum\limits_{v\in V(G)}\tr_G(v)$, where $W(G)$ is  the famous Wiener index of $G$. The \textit{Wiener complexity} $C_W(G)$ of a graph $G$ was introduced in~\cite{AAKS2014} (under the name Wiener dimension) as the number of different transmission of vertices in $G$:
$$C_W(G) = |\{ \tr_G(v):\ v\in V(G) \}|\,.$$
The Wiener complexity of graphs has been further investigated in~\cite{alizadeh-2019, AK2018, jamilet-2016, KJRMP2018+, XIIK2020+}.  Complexities of related invariants of interest in mathematical chemistry have also been investigated; the complexity of the connective eccentric index in~\cite{AK2016, CL2019}, the eccentric complexity in~\cite{ADX2017+}, and the complexity of the Szeged index in~\cite{AK2020}.

The \textit{transmission set} $\tr(G)$ of $G$ is the set of the transmissions of its vertices, that is, $\tr(G) = \{\tr_G(v):\ v\in V(G)\}$. A graph $G$ is \textit{transmission regular}~\cite{LDW2016} if all its vertices have the same transmission. In other words, transmission regular graphs are the graphs $G$ with $C_W(G) = 1 = |\tr(G)|$. On the other extreme, $G$ is \textit{transmission irregular}~\cite{AK2018} if all its vertices have pairwise different transmissions, that is, if $C_W(G)=n(G) = |\tr(G)|$.  We note in passing that very recently stepwise transmission irregular graphs were introduced in~\cite{DS2020} as the graphs in which the transmissions of any two of its adjacent vertices differ by exactly one.

Now, since almost no graphs are transmission irregular~\cite{AK2018}, it is of interest to search for families of transmission irregular graphs.  For this sake let $t\ge 3$, and  let  $k_1\, \ldots, k_t$ be positive integers. Then a {\em starlike tree} $T(k_1,\ldots,k_t)$ is a tree obtained by attaching to an isolated vertex $t$ pendant paths of lengths $k_1,\ldots,k_t$, respectively. These pendant paths will be called  {\em $k_i$-arms}. We may assume without loss of generality throughout the paper that $k_1\le \cdots \le k_t$. In~\cite{AK2018} it was proved that $T(1,k_2,k_3)$ is transmission irregular if and only if $k_3=k_2+1$ and $k_2\notin\left\{(t^2-1)/2,(t^2-2)/2\right\}$ for some $t\ge 3$.  Al-Yakoob and Stevanovi\'{c}~\cite{al-yakoob-2020} recently extended the latter result by characterizing the starlike trees $T(k_1,k_2,k_3)$ which are transmission irregular, their result will be restated in Theorem~\ref{thm:stevanovic}. In the meantime, Dobrynin constructed several families of transmission irregular graphs. In~\cite{dobrynin-2019} he presented an infinite family of $2$-connected transmission irregular graphs, in~\cite{dobrynin-2019-c} he followed with an infinite family of $3$-connected cubic transmission irregular graphs, while in~\cite{dobrynin-2019-b} he discovered an infinite family of transmission irregular trees of even order.

In the rest of this section we recall a few definitions needed and prove some preliminary results. In the first main result of Section~\ref{sec:arithmetic} we prove that if a starlike tree $T=T(a,a+1,\ldots,a+k)$, $k\ge 2$, is of odd order, then $T$ is transmission irregular. In the second main result of the section we then prove that a starlike tree $T(1,2,\ldots,\ell)$, $\ell \ge 3$, is transmission irregular if and only if $\ell \notin \{r^2 + 1:\ r\ge 2\}$. Then, in Section~\ref{sec:additional}, we determine additional infinite families among the starlike trees (broken unit arithmetic starlike trees and extremal starlike trees) and bi-starlike trees.  In the subsequent section we turn out attention to unicyclic graph containing $C_3$. From the two results proved we select the one asserting that the line graph of $T(a,a+1,a+2)$, $a\ge 2$, is transmission irregular if and only if $a$ is even.

\subsection{Preliminaries}

If $k$ is a positive integer, then $[k]=\{1, \ldots, k\}$ and $[k]_0=\{0,1,2,\ldots,k\}$. The degree of a vertex $v$ of a graph $G$ is denoted by $\deg_G(v)$. A vertex in a tree $T$ of degree at least $3$ is called a \textit{branching vertex} in $T$. The line graph of a graph $G$ is denote by $L(G)$. For an edge $e=uv$ of a graph $G$, the number of vertices that are closer to $u$ than to $v$ is denoted by $n_u(e|G)$ or $n_u$ for short. Analogously, $n_v(e|G)$ or $n_v$ for short denotes the number of vertices closer to $v$ than to $u$ in $G$. If $A$ is a set of integers and $i\in {\mathbb  Z}$, then $A+i$ denotes the usual coset, that is, $A + i = \{a+i:\ a\in A\}$.

We will make use of the following easy result on the transmission.

\begin{lemma}\label{equal} {\rm (\cite{Bala})}
If $uv\in E(G)$, then $\tr(u) - \tr(v) = n_v-n_u$.
\end{lemma}

If $T$ is a tree, then $n_u+n_v = n(T)$ for any edge $uv\in E(T)$. Hence in every tree $T$ there is at most one edge $uv$ for which $n_u = n_v$ holds. Moreover, if such an edge exists, then $n(T)$ must be even. Combining this fact with Lemma~\ref{equal}, we have the following result.

\begin{proposition}\label{TI-no-eq}
If $T$ is a transmission irregular tree, then $T$ contains no edge $uv$ with $n_u = n_v$.
\end{proposition}

Using Lemma~\ref{equal} we also derive the following result.

\begin{proposition}\label{TI-no-ad}
If $T$ is a transmission irregular tree, then $T$ contains no two edges $e_1 = xy$ and $e_2 = uv$ with $|n_x - n_y| = |n_u - n_v| = 1$.
\end{proposition}

\proof
From Lemma~\ref{equal} we get $\tr(x) = \tr(y) + n_y - n_x$. Since $|n_x - n_y| = |n_u - n_v| = 1$, the edges $xy$ and $uv$ must be adjacent in $T$. We may thus assume without loss of generality that $e_2=yz$, where $z\ne x$. Then $|n_x - n_y| = |n_z - n_y| = 1$. Since $n_x + n_y =  n(T) = n_z + n_y$, we get that $n_x = n_z$. Hence, using Lemma~\ref{equal} again, we get $\tr(z) = \tr(y) + n_y - n_z = \tr(y) + n_y - n_x = \tr(x)$,  contradicting the assumption that $T$ is transmission irregular.  \qed

\begin{proposition}\label{Pro:pend}
Let $G$ be a connected graph with  $n(G)=n$ and $v\in V(G)$ of degree $\deg(v)\geq 3$. If $P=uv_1v_2\cdots v_{x-1}v$ is a pendant path with natural adjacency relation attaching at $v$, where $\deg(u)=1$ and $x<\frac{n}{2}$, then $\tr(v_{x-1})-\tr(v)=n-2x$.
\end{proposition}

\proof
By definition, we have $n_u=1$ and $n_{v_1}=n-1$, that is, $n_{v_1}-n_u=n-2$.  Similarly, $n_{v_2}-n_{v_1}=n-4$, $\ldots$, $n_{v_{x-1}}-n_{v_{x-2}}=n-2(x-1)$, and $n_v-n_{v_{x-1}}=n-2x>0$. By Lemma~\ref{equal}, we get $\tr(v_{x-1})-\tr(v)=n-2x$.
\qed

We conclude the preliminaries with the following necessary condition for transmission irregular starlike trees.

\begin{proposition}\label{prop:nece}
If $T(k_1,\ldots,k_t)$ is transmission irregular, then $k_t\leq \sum\limits_{i=1}^{t-1}k_i$.
\end{proposition}

\proof
Set $T=T(k_1,\ldots,k_t)$ and suppose on the contrary that $k_t > \sum_{i=1}^{t-1}k_i$. Since $n(T) = 1 + \sum_{i=1}^{t}k_i$, we get that $k_t>\frac{n(T)}{2}$. Let $P$ be the $k_t$-arm in $T(k_1,\ldots,k_t)$. Based on the parity of $n(T)$, we observe that there exists an edge $uv$ on $P$ with $n_u = n_v$ if $n(T)$ is even, or there are two adjacent edges $xy$ and $yz$ with $|n_x - n_y| = |n_y - n_z| = 1$. By Propositions~\ref{TI-no-eq} and~\ref{TI-no-ad}, $T$ is not transmission irregular. This contradiction completes the proof.
\qed

\section{Unit arithmetic starlike trees}
\label{sec:arithmetic}

As already mentioned in the introduction, transmission irregular trees $T(1,k_2,k_3)$ were characterized in~\cite{AK2018}, while in~\cite{al-yakoob-2020} the result was extended to all starlike trees $T(k_1,k_2,k_3)$. We now restate this appealing result to show that the problem is intricate, as well as to be applied later on. Note that its condition $k_3 \le k_1 + k_2$ is just the case $t=3$ of Proposition~\ref{prop:nece}.

\begin{theorem}
\label{thm:stevanovic}  {\rm \cite[Theorem 2]{al-yakoob-2020}}
$T(k_1,k_2,k_3)$ is transmission irregular if and only if $k_1 < k_2 < k_3$, $k_3 \le k_1 + k_2$, and the triplet $(k_1,k_2,k_3)$ does not belong to the set
$$\bigcup_{1\le i<j} {\cal N}_{ij}^{xy} \cup \bigcup_{1\le j<k} {\cal N}_{jk}^{yz} \cup \bigcup_{1\le i<k} {\cal N}_{ik}^{xz}\,,$$
where
\begin{align*}
{\cal N}_{ij}^{xy} & = \Big\{ \Big(k_1, k_1+(j-i)\Big(1 + \frac{p}{\gcd(i+j,j-i)}\Big), \frac{p(i+j)}{\gcd(i+j,j-i)} \Big): \\
& i \le k_1,\ \frac{(k_1+j-i) \gcd(i+j,j-i)}{2i}\le p \Big\}\,, \\
{\cal N}_{jk}^{yz} & = \Big\{ \Big( \frac{p(j+k)}{\gcd(j+k,k-j)}, k_2, k_2 + (k-j) \Big(1+ \frac{p}{\gcd(j+k,k-j)}\Big) \Big) : \\
&  \max \Big(j, \frac{j+k}{\gcd(j+k,k-j)}\Big) \le k_2,\ 1\le p\le \frac{k_2\gcd(j+k,k-j)}{j+k} \Big\}\,, \\
{\cal N}_{ik}^{xz} & = \Big\{ \Big( k_1, \frac{p(i+k)} {\gcd(i+k,k-i)}, k_1 + (k-i) \Big( \frac{p} {\gcd(i+k,k-i)}\Big) \Big): \\
& i\le k_1,\ \frac{k_1 \gcd(i+k,k-i)}{i+k} \le p \le \frac{(k_1 + k - i) \gcd(i+k,k-i)}{2i} \Big\} \,.
\end{align*}
\end{theorem}

As pointed out by Al-Yokoob and Stevanovi\'{c}, the proof method used to prove Theorem~\ref{thm:stevanovic} could in principle be applied also to starlike trees with more than three arms. However, the number of sets of parameter values to be avoided grows quadratically with the number of branches, so possible formulations of such results (as well as their proofs) would be extremely long and consequently useless. Moreover, computational results (see~\cite[Table~1]{al-yakoob-2020}) indicate that the number of transmission irregular starlike trees rapidly decreases with the number of branches. Nevertheless we will construct in this section infinite families of transmission irregular starlike trees with an arbitrary number of arms.

We say that a starlike tree $T(k_1,\ldots,k_t)$ is \textit{arithmetic} if $k_{i+1}- k_i$ is a constant, or \textit{unit arithmetic} if $k_{i+1}-k_i=1$, for $i\in [t-1]$. The main result of this section reads as follows.

\begin{theorem}
\label{thm:Ti-a-k}
If $T$ is a unit arithmetic starlike tree of odd order, then $T$ is transmission irregular.
\end{theorem}

\proof
Let $T=T(a,a+1,\ldots,a+k)$, $k\ge 2$, and let $v$ be the vertex of $T$ with degree $k+1$. For $p\in [a+k]$ define the sets $B_p$ as follows:
$$B_p=\left\{
\begin{array}{ll}
\{ps+p(p-1)+2pi:\ i\in [k+1]\}; & p\in [a], \\
\\
\{ps+p(p-1)+2pi:\ i\in [(k+1)-(p-a)]\}; & p\in [a+k]\setminus [a]\,,
\end{array}
\right.$$ where $s=(k-1)(a+\frac{k}{2}-1)-2$.

\medskip\noindent
{\bf Claim A}: $\tr(T)\setminus\{\tr(v)\}=\bigcup\limits_{p=1}^{a+k}\Big(B_p+(\tr(v)+s+2)\Big)$. \\
By the structure of $T$ we see that $n(T)=1+a+(a+1) + \cdots + (a+k) = (k+1)(a+\frac{k}{2})+1$. Since the distance  between $v$ and the leaf on the longest arm is is $a+k$, we have $\tr(v^{\prime})-\tr(v)=(k+1)(a+\frac{k}{2})+1-2(a+k)=(k-1)(a+\frac{k}{2}-1)$ by Proposition~\ref{Pro:pend}, where $v^{\prime}$ is the neighbor of $v$ lying on the longest arm of $T$.  From Lemma \ref{equal}, the transmission of the vertex $w$ is $ps+p(p-1)+2pi$ if $w$ is on the $(a+k+1-i)$-arm of $T$ with $d(w,v)=p$ and $i\in [k+1]$  for $p\in [a]$ or $i\in [(k+1)-(p-a)]$ for $p\in [a+k]\setminus [a]$ where $s=(k-1)(a+\frac{k}{2}-1)-2$. Claim A now follows from the definition of $B_p$.
\smallqed

\medskip
By the assumption, the order of $T$, $n(T)=(k+1)(a+\frac{k}{2})+1$, is odd.
Note that $\tr(u)\neq\tr(v)$ for any vertex $u\in V(T)\setminus\{v\}$. Then it suffices to prove that $$|\{\tr(u):u\in V(T)\setminus\{v\}\}|=n(T)-1.$$

Note that $|A+a|=|A|$ for any set $A$. By Claim~A and the definition of $B_p$, it suffices to prove that $|\bigcup\limits_{p=1}^{a+k}B_p|=n(T)-1$, that is, the sets $B_p$, $p\in [a+k]$, are pairwise disjoint. Recall that $s=(k-1)(a+\frac{k}{2}-1)-2$ and note that $s$ is odd. Then $B_p$ consists of increasingly odd numbers in terms of $i$ if $p$ is odd, or of increasingly even numbers in terms of $i$ if $p$ is even. Set $B^{(1)}=\{B_p:p\in [a+k] \mbox{~is odd}\}$ and $B^{(2)}=\{B_p:p\in [a+k] \mbox{~is even}\}$. Since $B^{(1)}\cap B^{(2)}=\emptyset$, it suffices to prove that $B_p\cap B_{p+2t}=\emptyset$ for any subset $\{p,p+2t\}\subseteq [a+k]$. Next we calculate the value of  $\min B_{p+2}-\max B_{p}$. If $\{p,p+2\}\subseteq [a]$, we have
\begin{eqnarray*}
\min B_{p+2}-\max B_{p}
&=&(p+2)s+(p+2)(p+1)+2(p+2)\\
&&-ps-p(p-1)-2p(k+1) \\
&=&2s+6p+6-2p(k+1) \\
&=&(k-1)(2a+k-2)+6p+2-2p(k+1)\\
&=&(k-1)(2a-2p+k-2)+2p+2 \\
&>&0 \,.
\end{eqnarray*} If $p\in [a]$ and $p+2\in [a+k]\setminus [a]$, similarly as above,  we can get $\min B_{p+2}-\max B_p>0$. While $\{p,p+2\}\subseteq [a+k]\setminus [a]$, we have \begin{eqnarray*}
\min B_{p+2}-\max B_{p}
&=&(p+2)s+(p+2)(p+1)+2(p+2)\\
&&-ps-p(p-1)-2p[(k+1)-(p-a)] \\
&=&2s+6p+6-2p[(k+1)-(p-a)] \\
&=&(k-1)(2a+k-2)+6p+2-2p(k+1)+2p(p-a)\\
&=&(k-1)\Big[k-2-2(p-a)\Big]-4p+6p+2+2p(p-a)\\
&=&(k-1)(k-2)+2(p-a)(p-k+1)+2p+2.
\end{eqnarray*} Set $x=(k-1)(k-2)+2(p-a)(p-k+1)+2p+2$. Note that $k\geq 3$ and $a<p\leq a+k-2$. If $a\geq k$ or $a<k\leq p$, then $x>0$ holds clearly. Now we consider the last case $a<p<k$. In this case, since $k-p\geq 1$, we have \begin{eqnarray*}
x&=&(k-1)(k-2)-2(k-p)(p-a)+4p-2a+2\\
&=&k^2-3k+2-2(kp-p^2-ak+ap)+4p-2a+2 \\
&=&(k-p)^2+p^2+2a(k-p)+4p-2a-3k+4\\
&\geq&2p(k-p)+2a(k-p)-3(k-p)+p-2a+4\\
&=&(2p+2a-3)(k-p)+p+4-2a\\
&\geq& 2p+2a-3+p+4-2a\\
&=&3p+1\\
&>&0.
\end{eqnarray*}
Thus the sets $B_p$, $p\in [a+k]$, are pairwise disjoint, completing the proof.
\qed

It can be routinely checked that $(k+1)(a+\frac{k}{2})+1$ is odd if and only if $k\equiv 3({\rm mod}~4)$ or $k+2a\equiv 0({\rm mod}~4)$.
 Therefore, we have the following consequence.

\begin{corollary}\label{co:T-a-k} If $k\equiv 3({\rm mod}~4)$ or $k+2a\equiv 0({\rm mod}~4)$, then $T(a,a+1,\ldots,a+k)$ is transmission irregular.\end{corollary}

From Corollary~\ref{co:T-a-k} we can obtain some special transmission irregularity starlike trees. For instance,  $T(a,a+1,a+2,\ldots,a+k)$  is transmission irregular when $k\equiv 2({\rm mod}~4)$ and $a$ is odd. This fact for $k=2$ and odd $a$ enlarges the set of transmission irregular starlike trees  included in \cite[Table~1]{AK2018}. Moreover, we also have the following characterization for this case.

\begin{corollary}\label{com-k=2}
$T(a,a+1,a+2)$ is transmission irregular if and only if $a$ is odd.
\end{corollary}

\proof By the above we only need to prove that $T(a,a+1,a+2)$ is not transmission irregular if $a$ is even. Assume that $a=2t$ with $t\geq 1$. Setting $a=2t$ and $k=2$ in $B_p$, we have
$$B_p=\left\{
\begin{array}{ll}
\{ps+p(p-1)+2pi:i\in [3]\}; & p\in [2t], \\
\\
\{ps+p(p-1)+2pi:i\in [2t+3-p]\}; & p\in \{2t+1,2t+2\}\,,
\end{array}
\right.$$
where $s=2t-2$. Thus $\min B_{t+1}=(t+1)(s+t+2)=\max B_{t}$, which implies that $T(a,a+1,a+2)$ is not transmission irregular for $a=2t$ with $t\geq 1$.
\qed

Corollary~\ref{com-k=2} can also be deduced from Theorem~\ref{thm:stevanovic}. In the theorem, set $k_1 = a$, $k_2 = a+1$, and $k_3 = a+2$. Then it is easily seen that the sets ${\cal N}_{ij}^{xy}$ and $ {\cal N}_{jk}^{yz}$ are empty. For the set ${\cal N}_{ik}^{xz}$ we get that $k-i = 1$ and $p = \gcd(i+k,k-i) = 1$ must hold. For the second coordinate we have $i+k = a+1$, from which we get $2i = a$. Hence $(a,a+1,a+2)\in {\cal N}_{ik}^{xz}$ if and only if $a$ is even which, by Theorem~\ref{thm:stevanovic}, in turn implies that  $T(a,a+1,a+2)$ is transmission irregular if and only if $a$ is odd.

In the next result we provide a complete characterization of the transmission irregularity of arithmetic starlike trees with $k_1 = 1$.

\begin{theorem}\label{sl-from-one} Let $T=T(1,2,\ldots,\ell)$ with $\ell \ge 3$. Then $T$ is transmission irregular if and only if $\ell \notin \{r^2 + 1:\ r\ge 2\}$.
\end{theorem}

\proof Assume that $v$ is the vertex of maximum degree in $T$. Setting $a=1$ and $k=\ell-1$ in Claim A, we have $\tr(T)\setminus\{\tr(v)\}=\bigcup\limits_{p=1}^{\ell}\Big(B_p+(\tr(v)+s+2)\Big)$ where  $B_p=
\{ps+p(p-1)+2pi:i\in [\ell+1-p]\}$  with $s=(k-1)(a+\frac{k}{2}-1)-2=\frac{\ell(\ell-3)}{2}-1$ for $p\in [\ell]$. Clearly, $B_p$ is a set of increasing elements in terms of $i$. For convenience, we write $B_{p,i}=ps+p(p-1)+2pi$. Observe that $T$ is transmission irregular if and only if  $|\bigcup\limits_{p=1}^{\ell}B_p|=n(T)-1$, that is, the sets $B_p$, $p\in [\ell]$, are pairwise disjoint.  Let $B_p^{\prime}=B_p\setminus\{ps+p(p-1)+2p(\ell+1-p)\}$. Then we have \begin{eqnarray*}\min B_p-\max B_{p-1}^{\prime}&=&\frac{\ell(\ell-3)}{2}+p^2-(p-1)(p-3)-2(p-1)(\ell+1-p)\\
&=&\frac{\ell(\ell-3)}{2}+4p-3-2(p-1)(\ell+1-p) \\
&=&\frac{\ell(\ell-3)}{2}+2p^2-2\ell p+2\ell-1\\
&=&2\Big(p-\frac{\ell}{2}\Big)^2+\frac{\ell-2}{2}\\
&>&0,
\end{eqnarray*}
 that is, $\min B_p$ is larger than the second largest element in $B_{p-1}$ for any $p\in [\ell]\setminus \{1\}$. Note that \begin{eqnarray*}
\min B_p-\max B_{p-1}&=&\frac{\ell(\ell-3)}{2}+p^2-(p-1)(p-3)-2(p-1)(\ell+2-p)\\
&=&\frac{\ell(\ell-3)}{2}+2p^2-2(\ell+1)p+2\ell+1\\
&=&2\Big(p-\frac{\ell+1}{2}\Big)^2-\frac{\ell-1}{2}\\
&\geq&0
\end{eqnarray*} for any $p\in [1,\frac{\ell+1-\sqrt{\ell-1}}{2}]\cup[\frac{\ell+1+\sqrt{\ell-1}}{2},\ell]$. Moreover,  \begin{eqnarray*}\max B_{p-1}-\min B_p&=&\frac{\ell-1}{2}-2\Big(p-\frac{\ell+1}{2}\Big)^2\\
&\geq&0
\end{eqnarray*}  for any $p\in[\frac{\ell+1-\sqrt{\ell-1}}{2},\frac{\ell+1+\sqrt{\ell-1}}{2}]$. So $\max B_{p-1}\leq\min B_p\leq\min B_{p-1}^{\prime}$ for $p\in [\ell]\setminus \{1\}$. In view of Theorem \ref{thm:Ti-a-k}, we observe that $T$ is transmission irregular if and only if $n(T)$ is odd, or otherwise $\ell+1+\sqrt{\ell-1}$ is not even. Note that $n(T)=\frac{\ell(\ell+1)}{2}+1$ is odd if and only if $\ell \equiv j({\rm mod}~4)$ with $j\in \{0,3\}$.

We have thus proved that $T$ is \underline{not} transmission irregular if and only if $\ell \equiv j({\rm mod}~4)$ with $j\in \{1,2\}$ and $\ell+1+ \sqrt{\ell-1}$ is even. Suppose that  $\sqrt{\ell-1} = r \in {\mathbb Z}^+$. Then $\ell = r^2 + 1$ and hence $\ell + 1 + \sqrt{\ell-1} = (r^2+1) + 1 + r = r(r+1) + 2$, which is even. If $r=2k$, then $\ell = 4k^2 + 1$, so $\ell \equiv 1({\rm mod}~4)$. And if $r=2k+1$, then $\ell = 4k(k+1) + 2$, so $\ell \equiv 2({\rm mod}~4)$. We conclude that $T$ is not transmission irregular if and only if $\ell \in \{r^2 + 1:\ r\ge 2\}$.
\qed

To conclude the section we give a negative result by proving the a certain family of arithmetic starlike trees is not transmission irregular.

\begin{theorem}\label{Non-TI} If $\frac{2(a-3)}{3}\leq k\leq 2a+2$ and $k+2a\equiv 2({\rm mod}~ 4)$, then  $T(a,a+1,a+2,\ldots,a+k)$ is not transmission irregular.\end{theorem}
\proof Assume that $k+2a=4x+2$. Then $a+\frac{k}{2}-1=2x$. Since $\frac{2(a-3)}{3}\leq k\leq 2a+2$, we have $x\leq \min\{a,k+1\}$. Let $v$ be the vertex of maximum degree in $T(a,a+1,a+2,\ldots,a+k)$. As stated in the proof of Theorem \ref{thm:Ti-a-k}, we have
$\max B_x=xs+x(x-1)+2x(k+1)$ and $\min B_{x+1}=(x+1)s+x(x+1)+2(x+1)$ with $s=(k-1)(a+\frac{k}{2}-1)-2$. A straightforward calculation shows that $\max B_x=\min B_{x+1}$, which implies that $C_W(T)\leq n(T)-1$. \qed

\section{More (bi-)starlike transmission irregular trees}
\label{sec:additional}

\subsection{Broken unit arithmetic starlike trees}

A starlike tree is \textit{broken  unit arithmetic} if its arm length set is obtained from a  unit arithmetic sequence by removing  some consecutive elements of the sequence. If $a_1,\ldots, a_k$ is a unit arithmetic sequence in which  all elements from the open interval $(a_i,a_j)$ were removed, where $1\le i < j-1 \le k-1$, then the broken unit arithmetic starlike tree corresponding to this new sequence will be denoted by $T[a_1,a_i;a_j,a_k]$. Below we present the transmission irregularity of  a special class of broken arithmetic starlike trees $T[a,a+k-2;a+k,a+k+1]$.

\begin{theorem}
Let $T=T[a,a+k-2;a+k,a+k+1]$ with $k\geq 2$. If $k\equiv 3({\rm mod}~4)$ or $k+2a\equiv 0({\rm mod}~4)$, then $T$ is transmission irregular.
\end{theorem}

\proof Assume that $T_0=T(a,a+1,\ldots,a+k)$ with $u_0$ and $v_0$  being the pendant vertices of the $(a+k)$-arm and the $(a+k-1)$-arm  of $T_0$, respectively. Note that $u_0$ and $v_0$ are the diametrical vertices of $T_0$. Then $T$ can be obtained from $T_0$ by adding pendant vertices $u$ and $v$, and edges $uu_0$ and $vv_0$.  Let $w\in V(T)$ be the branching vertex of $T$  and let  $B_p$ and $s$  be defined as that in the proof of Theorem \ref{thm:Ti-a-k}. For $p\in [a+k]$ we define a new set $B_p^{\prime}$ which consists of the first two elements of $B_p$, that is, $B_p^{\prime}=\{ps+p(p-1)+2pi:\ i\in [2]\}$, and  define in addition the sets $B_p''$ by
$$B_p^{\prime\prime}=\left\{
\begin{array}{ll}
\{ps+p(p-1)+2(p+1)i:\ i\in [k+1]\}; & p\in [a], \\
\\
\{ps+p(p-1)+2(p+1)i:\ i\in [(k+1)-(p-a)]\}; & p\in [a+k-2]\setminus [a]\,,
\end{array}
\right.$$
with $s=(k-1)(a+\frac{k}{2}-1)-2$.
Let $B_p^*=B_p^{\prime}\cup B_p^{\prime\prime}$ for $p\in [a+k-2]$ and set $B_{a+k-1}^* = B_{a+k-1}^{\prime}$ for consistency.  The transmissions of vertices not on the diametrical path of $T$ form the set $\bigcup\limits_{p=1}^{a+k-2}(B_p^{\prime\prime}+\tr_T(w))$, those of vertices but $u$ and $w$ on the diametrical path of $T$ is just $\bigcup\limits_{p=1}^{a+k-1}(B_p^{\prime}+\tr_T(w))$. Let $B_{a+k}^*=\{\tr_T(u_0),\tr_{T}(v)\}$ and $B_{a+k+1}^*=\{\tr_T(u)\}$.  Note that $n(T)=n(T_0)+2$ is odd from the assumption.  Then
$$\tr(T)\setminus\{\tr_T(w)\}= \bigcup\limits_{p=1}^{a+k+1}B_p^*\,.$$
Moreover, $\tr_T(z)-\tr_T(w)=\tr_{T_0}(z)-\tr_{T_0}(w)$ for $z\in \{u_0,v_0\}$.  By Lemma \ref{equal}, we have  $\tr_T(u)=\tr_T(u_0)+n(T)-2$, and $\tr_T(v)=\tr_T(v_0)+n(T)-2$. From the proof of Theorem \ref{thm:Ti-a-k}, we know that $\tr_{T_0}(u_0)$ and $\tr_{T_0}(v_0)$ are maximum and second maximum transmissions in $\tr(T_0)$, so are $\tr_T(u)$ and $\tr_T(v)$ in $\tr(T)$. Now we only need to prove that   the sets from the union $\bigcup\limits_{p=1}^{a+k+1}B_p^*$  are pairwise disjoint. By a similar reasoning as that in the proof of Theorem \ref{thm:Ti-a-k}, it suffices to prove that  $\min B_{p+2}^*-\max B_{p}^*>0$ for any $\{p,p+2\}\subseteq [a+k]$. Let $x=\min B_{p+2}^*-\max B_{p}^*$. For any $\{p,p+2\}\subseteq [a]$, we have \begin{eqnarray*}
x&=&(k-1)(2a-2p+k-2)+2p+2-2(k+1) \\
&=&(k-1)[2a-2(p+2)+k]+2p-2\\
&\geq&k(k-1)+2p-2>0 \,.
\end{eqnarray*}
 Similarly, we have $x>0$ if $p\in [a]$ and $p+2\in [a+k]\setminus [a]$. For any $\{p,p+2\}\subseteq [a+k]\setminus [a]$, we have \begin{eqnarray*}
x&=&(k-1)(k-2)+2(p-a)(p-k+1)+2p+2-2[(k+1)-(p-a)] \\
&=&(k-1)(k-4)+2(p-a)(p-k+2)+2p-2.
\end{eqnarray*} If $a\geq k$ or $a<k\leq p$, then $x>0$ holds.  If $a<p<k$, we have  \begin{eqnarray*}
x&=&(k-p)^2+p^2+(2a-5)(k-p)+p-4a+2 \\
&\geq&(2p+2a-5)(k-p)+p-4a\\
&\geq&3p-2a-3\\
&\geq&3a+3-2a-3\\
&=&a>0,
\end{eqnarray*} completing the proof. \qed

\subsection{Extremal starlike trees}

In view of Proposition~\ref{prop:nece} we say that a starlike tree $T(k_1,\ldots,k_t)$ is {\em extremal} if $k_t=\sum_{i=1}^{t-1}k_i$ holds. In this section we construct some transmission irregular extremal starlike trees. Similarly as in the proof of Theorem~\ref{thm:Ti-a-k}, for positive integers $a$ and $k$ set where $h=(k-1)(2a+k)+2a-1$ and define the sets $D_p$, $p\in [a+k]$, as follows:
$$D_p=\left\{
\begin{array}{ll}
\{ph+p(p-1)+2pi:i\in [k+1]\}; & p\in [a], \\
\\
\{ph+p(p-1)+2pi:i\in [(k+1)-(p-a)]\}; & p\in [a+k]\setminus [a].
\end{array}
\right.$$
Mimicking the proof of Theorem \ref{thm:Ti-a-k} we can prove the following result, hence its proof is omitted.

\begin{lemma}\label{para}
Let $a$ and $k$ be positive integers, and let $D_p$ and $h$ be defined as above. Then the sets $D_p$, $p\in [a+k]$, are pairwise disjoint.
\end{lemma}

With Lemma~\ref{para} in hand we can find the announced transmission irregular extremal starlike trees.

\begin{theorem}\label{con-LDBS}
Let $T=T(a,a+1,\ldots,a+k,(a+\frac{k}{2})(k+1))$, and let $D=\bigcup\limits_{p=1}^{a+k}D_p$, where the sets $D_p$ are defined as above. If for every $d\in D$, the number $d$ is not a square number from the interval $\Big[k(2a+k-1)+1,(a+\frac{k}{2})^2(k+1)^2\Big]$, then $T$ is transmission irregular.
\end{theorem}

\proof  Note that $n(T)=(2a+k)(k+1)+1$. Assume that $v$ is the vertex with maximum degree in $T$ and $\tr(v)=y$. By Proposition \ref{Pro:pend}, we have $\tr(v_1)=y+1$ where $v_1$ lies on the longest pendant path in $T$ with $d_T(v,v_1)=1$. By Lemma \ref{equal}, the set of transmissions is just $\{1,4,9,\ldots,(a+\frac{k}{2})^2(k+1)^2\}+y$ of vertices on the longest pendant in $T$. Note that $h=(k-1)(2a+k)+2a-1$. From the structure of $T$, we observe that the set of transmissions is just $D+y$ of vertices in $T$ different from $v$ and not lying the longest arm. From the assumption with Lemma \ref{para}, our result follows.  \qed

Taking $k=1$ in Theorem \ref{con-LDBS}, we have  the following result.

\begin{corollary}\label{three-LDBS}
Let $T=T(a,a+1,2a+1)$. If $d$ is not a square number in the interval $\Big[2a+1,(2a+1)^2\Big]$ for any $d$ in $\{p(2a-1)+p^2:\ p\in [a+1]\}$ or $\{p(2a-1)+p^2+2p:\ p\in [a]\}$, then $T$ is transmission irregular.
\end{corollary}

\subsection{Bi-starlike trees}

A tree $T$ is a \textit{bi-starlike tree} if $T$ contains exactly two vertices of degrees at least $3$. The length of the induced path connecting these two vertices of degrees at least $3$ is called the \textit{shoulder width} of a bi-starlike tree and the path connecting these two vertices of degrees at least $3$ is called the \textit{shoulder path} in this bi-starlike tree. Denote by $BT^{(k_i)}(k_1,k_2,\ldots,k_t)$ a bi-starlike tree obtained from two copies, say $T_1$ and $T_2$, of starlike tree $T(k_1,k_2,\ldots,k_t)$ by identifying the branching vertex of $T_1$ with the leaf on the $k_i$-arm of $T_2$. Note that the shoulder width of $BT^{(k_i)}(k_1,k_2,\ldots,k_t)$ is $k_i$.

\begin{remark}\label{prime-bs} It can be routinely checked that $BT^{(2a+1)}(a,a+1,2a+1)$ is transmission irregular for $a\in \{2,3,5\}$ but not transmission irregular if $a=6$.\end{remark}

From Remark \ref{prime-bs}, it seems a bit difficult to construct transmission irregular bi-starlike trees with long shoulder widths. But the case is different when the shoulder width is short.
Although there are some unit arithmetic starlike trees which are not transmission irregular, we can construct transmission irregular bi-starlike trees with shoulder width $1$ using unit arithmetic starlike trees regardless of whether they are transmission irregular or not. Denote by $BS^*(a,a+1,a+2,\ldots,a+k)$ a bi-starlike tree obtained by connecting two vertices of degree $k+1$ of two copies of $T(a,a+1,a+2,\ldots,a+k)$ and attaching a pendant vertex to one vertex of degree $k+2$.
\begin{theorem}\label{wid=1} Let $T^*=BS^*(a,a+1,a+2,\ldots,a+k)$ with $a>1$. Then $T^*$ is transmission irregular.\end{theorem}
\proof From the structure of $T^*$, we have $n(T^*)=(k+1)(2a+k)+3$. Assume that $v,v^{\prime}\in V(T^*)$ with $\deg_{T^*}(v)=k+3$, $\tr(v)=x$ and   $\deg_{T^*}(v^{\prime})=k+2$. Then $vv^{\prime}\in E(T)$. By Lemma \ref{equal}, we have $\tr(v^{\prime})=x+1$. Let $T^*-vv^{\prime}=T\cup T^{\prime}$ where $v\in V(T)$ and $v^{\prime}\in V(T^{\prime})$. Setting $t={k\choose 2}+ak+1$, we have $n(T^*)=2t+1+2(a+k)$. Let $v_0$ be the leaf adjacent to $v$ in $T$ of $T^*$. Then $\tr(v_0)=x+2t+2(a+k)-1$. Now we define a set $A_i$ as follows:
 $$A_i=\left\{
\begin{array}{ll}
\{2it+i^2+2ij:j\in [k]_0\}; & i\in [a], \\
\\
\{2it+i^2+2ij:j\in [k+a-i]_0\}; & p\in [a+k]\setminus [a].
\end{array}
\right.$$

By Proposition \ref{Pro:pend}, the transmissions of vertices in $T$ adjacent to $v$ form the set $\Big(\{2t+2(a+k)-1\}\cup A_1\Big)+x$ and the transmissions of vertices in $T$ with distance $i$ to $v$ form the set $A_i+x$ for any $i\in [a+k]\setminus\{1\}$.  Moreover, $\tr(u^{\prime})=\tr(u)+1$ for any corresponding vertex $u^{\prime}$ in $T^{\prime}$ to $u$ in $T$ of $T^*$. Note that $|D|=|D+d|$ for any set $D$ and any number $d$. Then it suffices to prove that the sets  $A_1^*=A_1\cup\{2t+2(a+k)-1\}$ and $A_i$, $i\in \{2,\ldots a+k\}$, are pairwise disjoint.

 Since $a>1$, $A_1^*$ is pairwise disjoint. Moreover, by the definition of $A_i$, we have $\min\limits_{i\in [a+k]\setminus\{1\}}\min A_i=4t+4>2t+2(a+k)-1$ for $a>1$. Then we only need to prove that  $A_i$, $i\in \{2,\ldots a+k\}$, are pairwise disjoint. Note that $A_i$ consists of increasing odd numbers in terms of $i$  if $i$ is odd and vice versa. For any $\{i,i+2\}\subseteq[a]$, we have \begin{eqnarray*}
\min A_{i+2}-\max A_{i}&=&2(i+2)t+(i+2)^2-2it-i^2-2ik\\
&=&4t+4i+4-2ik\\
&\geq&2k(2a+k-1)-2a(k-2)+4\\
&>&0.
\end{eqnarray*}
From the fact that  $\max A_i=2it+i^2+2i(k+a-i)<2it+i^2+2ik$ for $i\in [a+k]\setminus [a]$, our results follows immediately. \qed

\section{Cycle-containing graphs}
\label{sec:Cycle-containing}

In this section we will construct some cycle-containing graphs with transmission irregularity. Denote by $C_3(k_1;k_2,k_3;k_4,k_5)$ a graph obtained from a triangle $C_3$ by attaching at one vertex of $C_3$ a pendant path of length $k_1$, at another vertex of $C_3$ pendant paths of lengths $k_2$ and $k_3$, respectively, and at the third vertex  pendant paths of lengths $k_4$ and $k_5$, respectively.

\begin{proposition}\label{k-al-3}
Let $k\ge 3$ and $G=C_3(1;1,k;2,k)$. Let $A_0=\{k+9,2k+11,2k+14,3k+14,4k+16\}$, $A_1=A+1$, $A=A_0\cup A_1$, and $B=\{i^2:\ i\in [k+3]\setminus[2]\}$. If $A\cap B=\emptyset$, then $G$ is transmission irregular.
\end{proposition}

\proof
Let $w$ be the unique vertex of degree $3$ in $G$, and let $u$ and $v$ be the two vertices of degree $4$ in $G$. From the structure of $G$, there is a pendant vertex $w^{\prime}$ attached at $w$ and there exist a pendant vertex $u^{\prime}$ and a pendant path $P_u:=uu_1u_2\ldots u_{k-1}u_{k}$ attached at $u$, two pendant paths $P^{\prime}:=vv^{\prime}v^{\prime\prime}$ and $P_v:=vv_1v_2\ldots v_{k-1}v_{k}$ attached at $v$ in $G$. Note that $n(G)=2k+7$. By the structure of $G$, we have $\tr(w)=k^2+3k+10$, $\tr(u)=(k+1)^2+9$, and $\tr(v)=(k+1)^2+8$. By  Proposition \ref{Pro:pend} and Lemma \ref{equal}, we observe that $\tr(w^{\prime})=k^2+5k+15$, $\tr(u^{\prime})=(k+1)^2+2k+14$, $\tr(v^{\prime})=(k+1)^2+2k+11$, $\tr(v^{\prime\prime})=(k+1)^2+4k+16$, the set of vertices on $P_u$ including $u$ is $\{(k+1)^2+j^2:j\in [k+3]\setminus [2]\}$ and the set of transmissions of vertices on $P_v$ including $v$ is $\{k^2+2k+j^2:j\in [k+3]\setminus [2]\}$. Therefore, we have
$$\tr(G)=D\bigcup (B+(k+1)^2)\bigcup (B+(k^2+2k)),$$
where $D=\{k+9,2k+11,2k+14,3k+14,4k+16\}+(k+1)^2$. Thus our result follows from the assumption.
\qed

Next we give a complete characterization of transmission irregularity of the line graph $L(T)$ of $T=T(a,a+1,a+2)$.

\begin{theorem}\label{cl-thr-star} Let $T=T(a,a+1,a+2)$ with $a\geq 2$. Then $L(T)$ is transmission irregular if and only if $a$ is even. \end{theorem}

\proof Note that $L(T)=C_3(a-1,a,a+1)$  of order $3a+3$. Assume that three vertices of degree $3$ in $L(T)$ are $u$, $v$, and $w$ at which the attached pendant paths are of lengths $a-1$, $a$, and $a+1$, respectively. Then $\tr(u)=\frac{a(3a+7)}{2}+4$, $\tr(v)=\frac{a(3a+7)}{2}+3$, and $\tr(w)=\frac{a(3a+7)}{2}+2$. From Proposition \ref{Pro:pend} and Lemma \ref{equal}, the transmissions of the vertices on the $(a-1)$-, $a$-, and $(a+1)$-arms not including $u$, $v$, and $w$ are $A_u+\frac{a(3a+7)}{2}$, $A_v+\frac{a(3a+7)}{2}$ and $A_w+\frac{a(3a+7)}{2}$, where $A_u=\{pa+(p+2)^2:p\in [a-1]\}$, $A_v=\{pa+(p+1)^2+2:p\in [a]\}$ and $A_w=\{pa+p^2+2:p\in [a+1]\}$. Let $A=\{2,3,4\}\cup A_u\cup A_v\cup A_w$. Then $\tr(L(T))=A+\frac{a(3a+7)}{2}$. Thus $L(T)$ is transmission irregular if and only if the sets $A_u$, $A_v$, and $A_w$ are pairwise disjoint. Since each of $A_u$, $A_v$, and $A_w$ consists of increasing numbers in terms of $p$, the only possible equal numbers in these three sets can happen if we have the following equality:
$$pa+(p+2)^2=(p+1)a+(p+1)^2+2,$$
which implies $2p+1=a$. We conclude that the sets $A_u$, $A_v$, and $A_w$ are pairwise disjoint if and only if $a$ is even.  \qed

From Theorem \ref{cl-thr-star}, we conclude that $T(a,a+1,a+2)$ is transmission irregular if and only if $LT(a,a+1,a+2)$ is not transmission irregular.  This interesting fact leads to the following problem.

\begin{problem}\label{TIS-L}
 Investigate the correlation between the transmission irregularity of (starlike) trees with that of their line graphs.
\end{problem}

\section*{Acknowledgements}
Kexiang Xu is supported by supported by NNSF of China (grant No.\ 11671202, and the China-Slovene bilateral grant 12-9). Sandi Klav\v{z}ar acknowledges the financial support from the Slovenian Research Agency (research core funding P1-0297, projects J1-9109, J1-1693, N1-0095, and the bilateral grant BI-CN-18-20-008).



\begin{thebibliography}{99}

\bibitem{AAKS2014}
  Y.~Alizadeh, V.~Andova, S.~Klav\v{z}ar, R.~\v{S}krekovski,
  Wiener dimension: Fundamental properties and (5,0)-nanotubical fullerenes,
  MATCH Commun.\ Math.\ Comput.\ Chem.\ 72 (2014) 279--294.

\bibitem{ADX2017+}
  Y.~Alizadeh, T.~Do\v{s}li\'c, K.~Xu,
  On the eccentric complexity of graphs,
  Bull.\ Malays.\ Math.\ Sci.\ Soc.\ 42 (2019) 1607--1623.

\bibitem{alizadeh-2019}
  Y.~Alizadeh, E.~Estaji, S.~Klav\v{z}ar, M.~Petkov\v{s}ek,
  Metric properties of generalized Sierpi\'nski graphs over stars,
  Discrete Appl.\ Math.\ 266 (2019) 48--55.

\bibitem{AK2016}
  Y.~Alizadeh, S.~Klav\v{z}ar,
  Complexity of topological indices: The case of connective eccentric index,
  MATCH Commun.\ Math.\ Comput.\ Chem.\ 76 (2016) 659--667.

\bibitem{AK2018}
  Y.~Alizadeh, S.~Klav\v{z}ar,
  On graphs whose Wiener complexity equals their order and on Wiener index of asymmetric graphs,
  Appl.\ Math.\ Comput.\ 328 (2018) 113--118.

\bibitem{AK2020}
  Y.~Alizadeh, S.~Klav\v{z}ar,
   Complexity of the Szeged index, edge orbits, and some nanotubical fullerenes,
   Hacet.\ J.\ Math.\ Stat.\ 49 (2020) 87--95.

\bibitem{al-yakoob-2020}
  S.~Al-Yakoob, D.~Stevanovi\'{c},
  On transmission irregular starlike trees,
  Appl.\ Math.\ Comput., to appear.

\bibitem{Bala}
  K.~Balakrishnan, M.~Changat, I.~Peterin, S.~\v{S}pacapan, P.~\v{S}parl, A.~R.~Subhamathi,
  Strongly distance-balanced graphs and graph products,
  European J.\ Combin.\  30 (2009) 1048--1053.

\bibitem{CL2019}
  X.~Chen, H.~Lian,
  Solution to a problem on the complexity of connective eccentric index of graphs,
  MATCH Commun.\ Math.\ Comput.\ Chem.\ 82 (2019) 133--138.

\bibitem{dobrynin-2019}
  A.~A.~Dobrynin,
  Infinite family of 2-connected transmission irregular graphs,
  Appl.\ Math.\ Comput.\ 340 (2019) 1--4.

\bibitem{dobrynin-2019-b}
  A.~A.~Dobrynin,
  Infinite family of transmission irregular trees of even order,
  Discrete Math.\ 342 (2019) 74--77.

\bibitem{dobrynin-2019-c}
  A.~A.~Dobrynin,
  Infinite family of $3$-connected cubic transmission irregular graphs,
  Discrete Appl.\ Math.\ 257 (2019) 151--157.

\bibitem{DS2020}
  A.~A.~Dobrynin, R.~Sharafdini,
 Stepwise transmission irregular graphs,
 Appl.\ Math.\ Comput.\  371 (2020)  paper no.\ 124949.

\bibitem{jamilet-2016}
  D.~A.~Jemilet, I.~Rajasingh,
  Wiener dimension of spiders, $k$-ary trees and binomial trees,
  Int.\ J.\ Pure Appl.\ Math.\ 109 (2016) 143--149.

\bibitem{KJRMP2018+}
  S.~Klav\v{z}ar, D.~A.~Jemilet,  I.~Rajasingh, P.~Manuel, N.~Parthiban,
  General transmission lemma and Wiener complexity of triangular grids,
  Appl.\ Math.\ Comput.\ 338 (2018) 115--122.

\bibitem{LDW2016}
  H.~Lin, K.~Ch.~Das, B.~Wu,
  Remoteness and distance eigenvalues of a graph,
  Discrete Appl.\ Math.\ 215 (2016) 218--224.

\bibitem{XIIK2020+}
  K.~Xu, A.~Ili\'{c}, V.~Ir\v{s}i\v{c}, S.~Klav\v{z}ar, H.~Li,
  Comparing Wiener complexity with eccentric complexity,
  submitted.

\end{thebibliography}
\end{document}